\documentclass[12pt]{article} 
\usepackage{latexsym} 
\usepackage{amssymb} 
\usepackage{amsmath} 

\begin{document}


\newtheorem{theorem}{Theorem} 
\newtheorem{problem}{Problem} 
\newtheorem{definition}{Definition} 
\newtheorem{lemma}{Lemma} 
\newtheorem{proposition}{Proposition} 
\newtheorem{corollary}{Corollary} 
\newtheorem{example}{Example} 
\newtheorem{conjecture}{Conjecture} 
\newtheorem{algorithm}{Algorithm} 
\newtheorem{exercise}{Exercise} 
\newtheorem{xample}{Example} 
\newtheorem{remarkk}{Remark} 
 
\newcommand{\be}{\begin{equation}} 
\newcommand{\ee}{\end{equation}} 
\newcommand{\bea}{\begin{eqnarray}} 
\newcommand{\eea}{\end{eqnarray}} 
\newcommand{\beq}[1]{\begin{equation}\label{#1}} 
\newcommand{\eeq}{\end{equation}} 
\newcommand{\beqn}[1]{\begin{eqnarray}\label{#1}} 
\newcommand{\eeqn}{\end{eqnarray}} 
\newcommand{\beaa}{\begin{eqnarray*}} 
\newcommand{\eeaa}{\end{eqnarray*}} 
\newcommand{\req}[1]{(\ref{#1})} 
 
\newcommand{\lip}{\langle} 
\newcommand{\rip}{\rangle} 
\newcommand{\uu}{\underline} 
\newcommand{\oo}{\overline} 
\newcommand{\La}{\Lambda} 
\newcommand{\la}{\lambda} 
\newcommand{\eps}{\varepsilon} 
\newcommand{\om}{\omega} 
\newcommand{\Om}{\Omega} 
\newcommand{\ga}{\gamma} 
\newcommand{\ka}{\kappa}
\newcommand{\rrr}{{\Bigr)}} 
\newcommand{\qqq}{{\Bigl\|}} 
 
\newcommand{\dint}{\displaystyle\int} 
\newcommand{\dsum}{\displaystyle\sum} 
\newcommand{\dfr}{\displaystyle\frac} 
\newcommand{\bige}{\mbox{\Large\it e}} 
\newcommand{\integers}{{\Bbb Z}} 
\newcommand{\rationals}{{\Bbb Q}} 
\newcommand{\reals}{{\rm I\!R}} 
\newcommand{\realsd}{\reals^d} 
\newcommand{\realsn}{\reals^n} 
\newcommand{\NN}{{\rm I\!N}} 
\newcommand{\DD}{{\rm I\!D}} 
\newcommand{\degree}{{\scriptscriptstyle \circ }} 
\newcommand{\dfn}{\stackrel{\triangle}{=}} 
\def\complex{\mathop{\raise .45ex\hbox{${\bf\scriptstyle{|}}$} 
     \kern -0.40em {\rm \textstyle{C}}}\nolimits} 
\def\hilbert{\mathop{\raise .21ex\hbox{$\bigcirc$}}\kern -1.005em {\rm\textstyle{H}}} 
\newcommand{\RAISE}{{\:\raisebox{.6ex}{$\scriptstyle{>}$}\raisebox{-.3ex} 
           {$\scriptstyle{\!\!\!\!\!<}\:$}}} 
 
\newcommand{\hh}{{\:\raisebox{1.8ex}{$\scriptstyle{\degree}$}\raisebox{.0ex} 
           {$\textstyle{\!\!\!\! H}$}}} 

\newcommand{\OO}{\won} 
\newcommand{\calA}{{\cal A}} 
\newcommand{\calB}{{\cal B}} 
\newcommand{\calC}{{\cal C}} 
\newcommand{\calD}{{\cal D}} 
\newcommand{\calE}{{\cal E}} 
\newcommand{\calF}{{\cal F}} 
\newcommand{\calG}{{\cal G}} 
\newcommand{\calH}{{\cal H}} 
\newcommand{\calK}{{\cal K}} 
\newcommand{\calL}{{\cal L}} 
\newcommand{\calM}{{\cal M}} 
\newcommand{\calO}{{\cal O}} 
\newcommand{\calP}{{\cal P}} 
\newcommand{\calX}{{\cal X}} 
\newcommand{\calXX}{{\cal X\mbox{\raisebox{.3ex}{$\!\!\!\!\!-$}}}} 
\newcommand{\calXXX}{{\cal X\!\!\!\!\!-}} 
\newcommand{\gi}{{\raisebox{.0ex}{$\scriptscriptstyle{\cal X}$} 
\raisebox{.1ex} {$\scriptstyle{\!\!\!\!-}\:$}}} 
\newcommand{\intsim}{\int_0^1\!\!\!\!\!\!\!\!\!\sim} 
\newcommand{\intsimt}{\int_0^t\!\!\!\!\!\!\!\!\!\sim} 
\newcommand{\pp}{{\partial}} 
\newcommand{\al}{{\alpha}} 
\newcommand{\sB}{{\cal B}} 
\newcommand{\sL}{{\cal L}} 
\newcommand{\sF}{{\cal F}} 
\newcommand{\sE}{{\cal E}} 
\newcommand{\sX}{{\cal X}} 
\newcommand{\R}{{\rm I\!R}} 
\renewcommand{\L}{{\rm I\!L}} 
\newcommand{\vp}{\varphi} 
\newcommand{\N}{{\rm I\!N}} 
\def\ooo{\lip} 
\def\ccc{\rip} 
\newcommand{\ot}{\hat\otimes} 
\newcommand{\rP}{{\Bbb P}} 
\newcommand{\bfcdot}{{\mbox{\boldmath$\cdot$}}} 
 
\renewcommand{\varrho}{{\ell}} 
\newcommand{\dett}{{\textstyle{\det_2}}} 
\newcommand{\sign}{{\mbox{\rm sign}}} 
\newcommand{\TE}{{\rm TE}} 
\newcommand{\TA}{{\rm TA}} 
\newcommand{\E}{{\rm E\,}} 
\newcommand{\won}{{\mbox{\bf 1}}} 
\newcommand{\Lebn}{{\rm Leb}_n} 
\newcommand{\Prob}{{\rm Prob\,}} 
\newcommand{\sinc}{{\rm sinc\,}} 
\newcommand{\ctg}{{\rm ctg\,}} 
\newcommand{\loc}{{\rm loc}} 
\newcommand{\trace}{{\,\,\rm trace\,\,}} 
\newcommand{\Dom}{{\rm Dom}} 
\newcommand{\ifff}{\mbox{\ if and only if\ }} 
\newcommand{\proof}{\noindent {\bf Proof:\ }} 
\newcommand{\remark}{\noindent {\bf Remark:\ }} 
\newcommand{\remarks}{\noindent {\bf Remarks:\ }} 
\newcommand{\note}{\noindent {\bf Note:\ }}

\newcommand{\boldx}{{\bf x}} 
\newcommand{\boldX}{{\bf X}} 
\newcommand{\boldy}{{\bf y}} 
\newcommand{\boldR}{{\bf R}} 
\newcommand{\uux}{\uu{x}} 
\newcommand{\uuY}{\uu{Y}} 
 
\newcommand{\limn}{\lim_{n \rightarrow \infty}} 
\newcommand{\limN}{\lim_{N \rightarrow \infty}} 
\newcommand{\limr}{\lim_{r \rightarrow \infty}} 
\newcommand{\limd}{\lim_{\delta \rightarrow \infty}} 
\newcommand{\limM}{\lim_{M \rightarrow \infty}} 
\newcommand{\limsupn}{\limsup_{n \rightarrow \infty}} 
 
\newcommand{\ra}{ \rightarrow }

\newcommand{\ARROW}[1] 
  {\begin{array}[t]{c}  \longrightarrow \\[-0.2cm] \textstyle{#1} \end{array} } 
 
\newcommand{\AR} 
 {\begin{array}[t]{c} 
  \longrightarrow \\[-0.3cm] 
  \scriptstyle {n\rightarrow \infty} 
  \end{array}} 
 
\newcommand{\pile}[2] 
  {\left( \begin{array}{c}  {#1}\\[-0.2cm] {#2} \end{array} \right) } 
 
\newcommand{\floor}[1]{\left\lfloor #1 \right\rfloor} 
 
\newcommand{\mmbox}[1]{\mbox{\scriptsize{#1}}} 
 
\newcommand{\ffrac}[2] 
  {\left( \frac{#1}{#2} \right)} 
 
\newcommand{\one}{\frac{1}{n}\:} 
\newcommand{\half}{\frac{1}{2}\:} 
 
\def\le{\leq} 
\def\ge{\geq} 
\def\lt{<} 
\def\gt{>} 
 
\def\squarebox#1{\hbox to #1{\hfill\vbox to #1{\vfill}}} 
\newcommand{\qed}{\hspace*{\fill} 
           \vbox{\hrule\hbox{\vrule\squarebox{.667em}\vrule}\hrule}\bigskip} 
 
\title{Sufficient conditions for the  invertibility of adapted perturbations of identity on the
  Wiener space}

\author{Ali  S\"uleyman  \"Ust\"unel and Moshe  Zakai} 
\date{ } 
\maketitle 
\noindent 
{\bf Abstract:}{\small{ Let $(W,H,\mu)$ be the classical  Wiener
    space. Assume that $U=I_W+u$ is an adapted perturbation of
    identity, i.e., $u:W\to H$ is adapted to the canonical filtration
    of $W$. We give some sufficient analytic  conditions on $u$
    which imply the invertibility of the map $U$. In particular it is
    shown that if $u\in \DD_{p,1}(H)$ is adapted and if
    $\exp(\frac{1}{2}\|\nabla u\|_2^2-\delta u)\in L^q(\mu)$, where
    $p^{-1}+q^{-1}=1$,  then
    $I_W+u$ is almost surely invertible. With the help of this result
    it is shown that if $\nabla u\in L^\infty(\mu,H\otimes H)$, then
    the Girsanov exponential of $u$ times the Wiener measure
    satisfies the logarithmic Sobolev inequality and this implies the
    invertibility of
    $U=I_W+u$.  As a consequence, if, there
    exists an integer $k\geq 1$ such that $\|\nabla^k
    u\|_{H^{\otimes(k+1)}}\in L^\infty(\mu)$, then $I_W+u$ is again
    almost surely invertible under the almost sure continuity
    hypothesis of $t\to\nabla^i \dot{u}_t$ for $i\leq k-1$. }}\\

\section{Introduction}
This paper is devoted to the search of sufficient conditions for the
invertibility of a certain class of mappings on the Wiener
space. This class consists of the mappings of the form of perturbation
of identity, where the perturbation part is a mapping which is
absolutely continuous with respect to the Lebesgue measure and the
corresponding density is adapted and almost surely square integrable. 
To make the things more precise,
let $W=C_0([0,1])$ be the Banach space of continuous functions on 
$[0,1]$, with its Borel sigma field denoted by $\calF$. We denote by 
$H$ the Cameron-Martin space, namely the space of absolutely 
continuous functions on $[0,1]$ with square integrable Lebesgue 
density:  
$$ 
H=\left\{h\in 
W:\,h(t)=\int_0^t\dot{h}(s)ds,\,|h|_H^2=\int_0^1|\dot{h}(s)|^2ds<\infty\right\}\,. 
$$ 
$\mu$ denotes the classical Wiener measure on $(W,\calF)$, 
$(\calF_t,t\in [0,1])$ is the filtration generated by the paths of the 
Wiener process $(t,w)\to W_t(w)$, where $W_t(w)$ is defined as $w(t)$ 
for $w\in W$ and $t\in [0,1]$.  
Assume that $u:W\to H$ is a measurable mapping, define  $U:W\to W$ as 
$$
U= I_W + u\,.
$$
 $U$  can be represented  as
\begin{equation} 
\label{eq:1.1} 
U_t(w) =W_t(w) + \int_0^t \dot{u}_s(w) ds,
\end{equation} 
using the isometry between $H$ and $L^2([0,1])$. We assume that  $\dot{u}$ is
adapted to the filtration $(\calF_t,\,t\in [0,1])$.  
For simplicity we consider in this paper the Banach space of continuous
functions on $[0,1]$, taking values in $\reals$; the results,
however, go over directly to the infinite dimensional case, including
the Wiener space corresponding to the cylindrical Wiener process  based
on a Hilbert space. 

\noindent
To illustrate a situation  where the addressed  problem comes up, consider the
question  of the absolute  continuity of the measure $U\mu$, i.e., the
image of $\mu$ under $U$  and the calculation of the  corresponding
Radon-Nikodym  derivative in case of absolute  continuity.  The
celebrated Girsanov theorem (cf.\cite{GIR, Orey}) yields the change of variables  
formula, i.e. setting 
$$ 
\rho_U(w) = \exp \left(- \int_0^1 \dot{u}_s dw_s - \half \int_0^1
  |\dot{u}_s|^2 ds \right) 
$$ 
and assuming that $E[\rho_U] = 1$, then, for smooth $f$, it holds true
that 
$$ 
E[f\circ U\, \rho_U] = E[f]\,. 
$$ 
Hence the image measure  $U\mu$ is absolutely continuous with respect
to $\mu$. Let $\la$ 
be the corresponding  Radon-Nikodym derivative:
$$ 
E[f\circ U] = E[f\,\la],  
$$ 
then  under ``suitable conditions''
\begin{equation} 
\label{eq:1.2} 
\la = \frac{1}{\rho_U\circ U^{-1}}
\end{equation} 
where $U^{-1}$ is the inverse to $U$ (cf. e.g. Section~1.3 of 
\cite{BOOK}). Therefore  the invertibility  of $U$ plays a fundamental
role in the evaluation of the Radon-Nikodym derivative  $\la$. This
situation is particularly important if we want to write a probability
density as the Radon-Nikodym derivative of the image of the Wiener
measure under a mapping of the form of perturbation of identity;   we
refer the reader  to \cite{FUZ} for a quick introduction to this problem.

\vspace{0.5cm} 
\noindent

The second and somehow related question concerns the question of the
existence and uniqueness of strong solutions of the  stochastic
differential equations of  the following type: 
\begin{eqnarray}
\label{sde}
dV_t&=&-\dot{u}_t\circ V+\,dW_t\\
 V_0&=&0\nonumber\,, 
\end{eqnarray}
where $\dot{u}$ is described above. Note here the fact that, though
adapted,  the drift coefficient $\dot{u}_t(w)$ may
depend on the whole history of the Brownian path  $(w(s),s\in [0,t])$. If one can show that the map
defined by $U=I_W+u$, where $u$ is the primitive of $\dot{u}$, has a left
inverse $V$,  then this  inverse map will be the unique solution
of the equation (\ref{sde}). In fact, under the hypothesis
$E[\rho_U]=1$ and  $\dot{u}\in L^2(dt\times d\mu)$,  we prove in Theorem
\ref{adaptedness-thm} that,  if $U=I_W+u$ has a left inverse
$V$, then the image of $\mu$ under $V$ is equivalent to $\mu$,
$V$ is also right inverse and it is of the form $V=I_W+v$ with  $v:W\to H$,
$\dot{v}$ is adapted and finally that $V$ is the unique strong
solution of the equation \ref{sde}. This result is quite useful and
seems to be new.

As it is shown by the  well-known counter example
given by Tsirelson (cf. \cite{I-W}, p. 181), the usual hypothesis of
integrability on $u$ does not imply the existence of strong solutions,
hence the invertibility of $U$ either. 
A well-known condition  for the  existence of such  an  inverse  is
the case where the drift coefficient is Lipschitz continuous in the
Cameron-Martin space direction. Namely 
$$
\sup_{s\leq t}|\dot{u}_s(w+h)-\dot{u}_s(w+k)|\leq K\,\sup_{s\leq
  t}|h(s)-k(s)|\,,
$$
$\mu$-a.s., for any $h,k\in H$, $t\in [0,1]$,  where $K$ is a
constant. In this case, using the usual fixed point techniques, one can
prove that the stochastic differential equation
has a unique adapted solution. Then it is clear that $V\mu$ is
equivalent  to $\mu$ and that $U\circ V=V\circ
U=I_W$ almost surely. In the sequel, between other things  we shall
also  surpass this frame.


\vspace{0.5cm}
\noindent
Let us summarize the contents of the paper:
the basic notions of functional analysis on the Wiener space and   the
stochastic calculus of variations (the Malliavin  
calculus) are reviewed in Section \ref{preliminaries}.  Section
\ref{suff-cond-1}   presents basic   results 
on the  invertibility of $U=I+u$ with  $u:W\to H$ adapted\footnote{For
  practical 
  reasons, we call $u$ adapted whenever its Lebesgue density $\dot{u}$, called
  sometimes the drift,  is adapted.}.  
The main results are obtained by the regularization
of the drift with the Ornstein-Uhlenbeck semigroup. In fact, let
$(P_\tau,\,\tau\geq 0)$ denote the Ornstein-Uhlenbeck 
semigroup (cf. the formula (\ref{O-U})) and set $e^{-\tau} P_\tau u=
u_\tau$.  It is shown that under  
reasonable  assumptions on $u$, the map  $U_\tau=I+u_\tau$ is
invertible. Its inverse is of the 
form $V_\tau=I_W+v_\tau$, where $v_\tau$ is $H$-valued and
adapted. Besides the following  identities are satisfied almost everywhere:
\begin{equation}
\label{inv-1}
v_\tau=-u_\tau\circ V_\tau
\end{equation}
and
\begin{equation}
\label{inv-2}
u_\tau=-v_\tau\circ U_\tau\,.
\end{equation}
If we can show that $v=\lim_{\tau\to 0}v_\tau$  exists in probability
and if it satisfies the relations  (\ref{inv-1}) and (\ref{inv-2})
where $u_\tau,\,U_\tau$ and $v_\tau,\,V_\tau$ are replaced
respectively by $u,\,U$ and $v,\,V$, 
then the  invertibility of $I+u$ follows. This program, which is
realized in Section \ref{suff-cond-1}, is not so
evident, in  fact
we need the Carleman inequality (cf. \cite{Carl,D-S}) to find useful
sufficient conditions to show the existence of this  limit. The basic
result of this ection proves that if $u\in \DD_{p,1}(H)$, $p>1$,  is adapted
and if $\exp(\frac{1}{2}\|\nabla u\|_2^2-\delta u)\in L^q(\mu)$, where
$p^{-1}+q^{-1}=1$, then
$U=I_W+u$ is almost surely invertible.  As an
application and a further tool also, we prove a logaritmic Sobolev
inequality for the measures of the type $d\nu=\rho_U\,d\mu$, where
$U=I_W+u$ and the Sobolev derivative of $u$ is essentially bounded as
a Hilbert-Schmidt operator. Using this inequality, we show that such a
$U$ is almost surely invertible. Although this is not the most
general sufficient  condition for the  invertibility  that we find,
as  hypothesis  it is strictly  weaker than
the Lipschitz  assumption , which is generally  used for 
the construction of the inverse mapping  via the stochastic differential
equations as illustrated with the formulae (\ref{sde}) (cf. Remark
\ref{remark-1} for the details).
  
Section \ref{extensions} extends these results using some localization
techniques. As a
corollary  we prove that, for any $k\geq 1$,  if the $k$-th order  Sobolev
derivative $\nabla^k u$ of $u$ is  essentially bounded as a Hilbert-Schmidt
tensor and if the Sobolev derivatives upto the order $k-1$ of the
process $t\to\dot{u}(t)$ 
are almost surely continuous, then $U$ is almost surely
invertible. For the case $k=1$ this continuity hypothesis is avoided
using the logarithmic Sobolev inequality as explained in Section
\ref{suff-cond-1}.

Finally we underline
the fact that the results of this paper can be extended to the
abstract Wiener spaces where the notion of adaptedness can be defined
with respect to any continuous resolution of identity of the
associated  Cameron-Martin space as indicated in \cite{FILT} or
sections 2.6 and 3.6 of \cite{BOOK}.

A preliminary version of these results have been announced in the note
\cite{CRAS}, however the contents of this paper are considerably
stronger and  more
general. In particular, using  a convex interpolation method we
succeed to diminish the Sobolev differentiability requirements about
the shift, to  prove a new logarithmic Sobolev inequality and several
extensions as explained in the last section of the paper.

\section{Preliminaries} 
\label{preliminaries}
Let $W=C_0([0,1])$ be the Banach space of continuous functions on
$[0,1]$, with its Borel sigma field denoted by $\calF$. We denote by
$H$ the Cameron-Martin space, namely the space of absolutely
continuous functions on $[0,1]$ with square integrable Lebesgue
density: 
$$
H=\left\{h\in
W:\,h(t)=\int_0^t\dot{h}(s)ds,\,|h|_H^2=\int_0^1|\dot{h}(s)|^2ds<\infty\right\}\,.
$$
$\mu$ denotes the classical Wiener measure on $(W,\calF)$,
$(\calF_t,t\in [0,1])$ is the filtration generated by the paths of the
Wiener process $(t,w)\to W_t(w)$, where $W_t(w)$ is defined as $w(t)$
for $w\in W$ and $t\in [0,1]$. We shall recall briefly  some well-known  functional
analytic tools on the Wiener space, we refer the reader to
\cite{Mal,F-P,ASU} or to \cite{ASU-1} for further details:
$(P_\tau,\tau\in \R_+)$ denotes the
semi-group of Ornstein-Uhlenbeck on $W$, defined as 
\begin{equation}
\label{O-U}
P_\tau f(w)=\int_Wf(e^{-\tau}w+\sqrt{1-e^{-2\tau}}y)\mu(dy)
\end{equation}
Let us recall that $P_\tau=e^{-\tau\calL}$, where $\calL$ is the
number operator. We denote by $\nabla$ the Sobolev derivative which is
the extension (with respect to the Wiener measure) of the Fr\'echet
derivative in the Cameron-Martin space
direction. The iterates of $\nabla$ are defined similarly. Note that,
if $f$ is real valued, then $\nabla f$ is a vector and if $u$ is an
$H$-valued map, then $\nabla u$ is a Hilbert-Schmidt  operator (on $H$)
valued map whenever defined. If $Z$ is a separable Hilbert space and
if $p>1,k\in \R$, we
denote by $\DD_{p,k}(Z)$ the $\mu$-equivalence classes of $Z$-valued
measurable mappings $\xi$, defined on $W$ such that
$(I+\calL)^{k/2}\xi$ belongs to $L^p(\mu,Z)$ and this set, equipped
with the norm
\begin{equation}
\label{eq-1}
\|\xi\|_{p,k}=\|(I+\calL)^{k/2}\xi\|_{L^p(\mu,Z)}
\end{equation}
becomes a Banach space. From the Meyer inequalities, we know that 
the norm defined by 
$$
\sum_{i=0}^k\|\nabla^i\xi\|_{L^p(\mu,Z\otimes H^{\otimes i})}\,,\,k\in
\NN\,,
$$
is equivalent to the norm $\|\xi\|_{p,k}$ defined by (\ref{eq-1}). 
We denote by $\delta$ the adjoint of $\nabla$ under $\mu$ and recall
that, whenever $u\in \DD_{p,0}(H)$ for some $p>1$ is
adapted\footnote{In the sequel the adapted elements of $\DD_{p,k}(H)$
  will be denoted by $\DD_{p,k}^a(H)$.}
, then
$\delta u$ is equal to the It\^o integral of the Lebesgue density of $u$:
$$
\delta u=\int_0^1\dot{u}_s dW_s\,.
$$
\noindent

\noindent
Let $X$ be a separable Hilbert space and let $f:W\to X$ be a
measurable map. We say that $f$ is an $H-C$-map if $f$ has a
modification (denoted again as $f$) such that the mapping $h\to
f(w+h)$ is continuous for $\mu$-almost all $w\in W$. Similarly, we say
that $f$ is $H-C^k$, $k\geq 1$ or that it is $H$-real analytic,  if $h\to f(w+h)$ is $k$-times
differentiable or real analytic $\mu$-almost surely. In the sequel, we
shall use the same notation for the $H$-derivative and for the Sobolev
derivative since the latter is the $L^p$-extension of the former. 
Note that the set 
$A=\{w\in W:\,h\to f(w+h)\in C^k(H)\}$ is $H$-invariant, i.e.,
$A+H\subset A$, hence $A^c$ has zero capacity as soon as $\mu(A)>0$ (cf. \cite{ASU-1}). The
following result is well-known (cf. \cite{BOOK}, Lemma 3.3.2):
\begin{lemma}
\label{regularity-lemma}
Assume that $f\in L^p(\mu,X)$, where $X$ is a separable Hilbert
space. Then, for any $\tau>0$, $P_\tau f$ has a modification $f_\tau$,
such that $h\to f_\tau(w+h)$ is almost surely analytic on $H$, in
other words $P_\tau f$ is $H$-analytic. In particular it is $H-C^\infty$.
\end{lemma}
 
\noindent
Another important result that we shall need is the following one
(cf. \cite{BOOK}, Theorem 3.5.3 where a more general case is treated
and  Theorem 4.4.1 in the $H-C^1$-case):
\begin{theorem}
\label{jacobi-thm}
Assume that $\xi:W\to H$ is an $H-C^1$-map and denote by $T$ the map $T=I_W+\xi$. Then
\begin{enumerate}
\item The set $T^{-1}\{w\}$ is countable $\mu$-almost surely. Let
  $N(w)$ be its cardinal.
\item For any $f,\,g\in C_b(W)$, we have  the change of variables
  formula:
$$
E[f\circ T\,g\,\rho_T]=E\left[f\, \sum_{y\in T^{-1}\{w\}}g(y)\right]\,, 
$$
in particular
$$
E[f\circ T\,\rho_T]=E[f\, N]\,, 
$$
where 
$$
\rho_T=\dett(I_H+\nabla\xi)\exp\left[-\delta\xi-\frac{1}{2}|\xi|_H^2\right]
$$
and $\dett(I_H+\nabla\xi)$ denotes the modified Carleman-Fredholm
determinant.
\end{enumerate}
\end{theorem}
\begin{remarkk}
If $A$ is a nuclear operator on a separable Hilbert space, then 
$\dett(I_H+A)$ is defined as 
\beaa
\dett(I_H+A)&=&\prod_{i=1}^\infty (1+\la_i)e^{-\la_i}\\
&=&\det(I_H+A)e^{-\trace A}\,,
\eeaa
where $(\la_i)$ denotes the spectrum of $A$ and each eigenvalue is
counted with respect to its multiplicity. Afterwards, one can show
that $A\to \dett(I_H+A)$ extends continuously (even analytically) to
the space of Hilbert-Schmidt operators, cf. \cite{D-S}. If $A$ is a
quasi-nilpotent operator, then by definition the spectrum of $A$
is equal to the singleton  $\{0\}$, hence, in this case we always have
$\dett(I_H+A)=1$. 
\end{remarkk}
\begin{remarkk}
It is well-known that (cf. \cite{Kus,BOOK}) given an  $H$-valued
$H-C^1$ map, there exists a measurable partition $(M_n,n\geq 1)$ of its set of
non-degeneracy $M$, i.e., the set on which the $\dett(I_H+\nabla u)$ is
non-zero, such that on each $M_n$, $I_W+u$ is equal to some invertible
mapping of the form of perturbation of identity. This result implies
that the notion of multiplicity $N$ is well-defined and it is equal to the
$$
N(w,M)=\sum_{y\in T^{-1}\{w\}\cap M}1_y\,.
$$
Besides, using the Sard Lemma on Wiener space  (cf. \cite{BOOK},
Proposition 4.4.1), one can show that 
$$
N(w,M)=N(w,W)
$$
almost surely. This result is  extended even  to $H-C^1_{\rm
  loc}$-maps (cf. Definition \ref{hc-1-defn})  as explained in Chapter
IV of \cite{BOOK}. 
\end{remarkk}

\begin{theorem}
\label{jacobi-bis}
 Let $\xi$ be as in Theorem \ref{jacobi-thm} with
 $\xi\in\DD_{p,1}(H)$ for some $p>1$. If the
 Lebesgue density of $\xi$, called drift and  denoted by
  $\dot{\xi}$ is adapted to the filtration of the canonical Wiener
  process, then $\rho_T$ reduces to the usual exponential martingale:
\beaa
\rho_T&=&\exp\left[-\delta\xi-\frac{1}{2}|\xi|_H^2\right]\\
&=&\exp\left[-\int_0^1\dot{\xi}_s\cdot
  dW_s-\frac{1}{2}\int_0^1|\dot{\xi}_s|^2ds\right]\,.
\eeaa
In this case we have always $N(w)\in \{0,1\}$ almost surely and $N=1$
a.s. if $E[\rho_T]=1$. 
\end{theorem}
\proof
The proof follows from the fact that $\delta \xi$ coincides with the It\^o
integral of $\dot{\xi}$ if the latter is adapted. In this case, we
always have, from the Fatou lemma $E[\rho_T]\leq 1$ and if
$E[\rho_T]=1$, then it follows from Theorem \ref{jacobi-thm} and from
the Girsanov theorem  that
$N=1$ almost surely.

\qed

\noindent
A simple, nevertheless important corollary of Theorems
\ref{jacobi-thm} and \ref{jacobi-bis}  is
\begin{corollary}
\label{cor-UZ}
Assume that $\xi$ is adapted and $H-C^1$. Assume moreover that
$E[\rho_T]=1$. Then, there exists a map $S$ of the form $S=I_W+\eta$
with  $\eta:W\to H$ adapted, such that
$$
\mu\left(\left\{w\in W:\,S\circ T(w)=T\circ S(w)=w\right\}\right)=1\,.
$$
In other words $T$ is almost surely invertible.
\end{corollary}
\proof
Let $\tilde{W}=\{w\in W:\,N(w)=1\}$, it follows from Theorem
\ref{jacobi-thm} and  Theorem \ref{jacobi-bis}  and from the
hypothesis about $\rho(-\delta \xi)$ that
$\mu(\tilde{W})=1$. Consequently, for any $w\in \tilde{W}$, there
exists a unique $S(w)\in W$ such that $T(S(w))=w$. Let us define $S$
on $\tilde{W}^c$ by $I_W$. Then, for any $A\in \calF$, we have 
$S^{-1}(A)=S^{-1}(A)\cap \tilde{W}$ $\mu$-almost surely. Moreover, by
the very definition of $S$, we have $S^{-1}(A)\cap \tilde{W}=T(A)\cap
\tilde{W}$. Since $\tilde{W}$ is a Borel set, to show the
measurability of $S$ with respect to the completion of $\calF$, it
suffices to show that $T(A)$ belongs to the same sigma-algebra and
this follows from Theorem 4.2.1 of \cite{BOOK}{\footnote{In fact $T(A)$  is  a
Souslin set as one can show by the help of the  measurable selection
theorem, hence  an element of the universal completion of
$\calF$, cf. \cite{Del-M}.}} and this settles the measurability of $S$. Moreover,
from Theorem \ref{jacobi-thm}, 
for any $f,\,g\in C_b(W)$, we have on the one hand 
\beaa
E[\rho(-\delta \xi)\,g]&=&E\left[\sum_{y\in
  T^{-1}\{w\}}g(y)\right]\\
&=&E[g\circ S]\,,
\eeaa
hence $S(\mu)$ is equivalent to $\mu$ and on the other hand
\beaa
E[f\circ S\circ T\,\rho(-\delta \xi)\,g]&=&E\left[f\circ S\,\sum_{y\in
  T^{-1}\{w\}}g(y)\right]\\
&=&E[f\circ S\,g\circ S]\\
&=&E[f\,g\,\rho(-\delta \xi)]\,.
\eeaa
Therefore $S\circ T=I_W$ almost surely. In particular $S$ is of the
form $I_W+\eta$ and $\eta$ is an  adapted and $H$-valued mapping.
\qed

\noindent
If $u\in \DD_{2,0}(H)$ is adapted and satisfies $E[\rho(-\delta
u)]=1$, but without $H-C^1$-hypothesis, in case it has a left inverse,
then this left  inverse is also a right inverse and it is  characterized
by the following Theorem:
\begin{theorem}
\label{adaptedness-thm}
Assume that $U=I_W+u$, $u:W\to H$, $u(t)=\int_0^t\dot{u}_sds$, for any
$t\in [0,1]$ and that $\dot{u}$ is adapted to the Brownian filtration
$(\calF_t,\,t\in [0,1])$. Assume further that $u\in L^2(\mu,H)$, 
$E[\rho(-\delta u)]=1$. Suppose that there exists some $V:W\to W$ such
that  $V\circ U=I_W$ a.s. Then 
\begin{enumerate}
\item $V\mu$ is equivalent to $\mu$ and $V$ is also a left inverse,
  i.e.,
$$
U\circ V=I_W
$$
$\mu$-almost surely. In other words $U$ is almost surely invertible
and its inverse is $V$. 
\item $V$ i s of the form of
a perturbation of identity, i.e.,  $V=I_W+v$ and  $v:W\to H$.
\item $\dot{v}$ is adapted to the filtration $(\calF_t,\,t\in [0,1])$.
\item In particular, the process $(V(t),\,t\in [0,1])$ is the unique
  strong solution of the stochastic differential equation (\ref{sde}).
\end{enumerate}
\end{theorem}
\proof
For any $f\in C_b(W)$, it follows  from the Girsanov theorem   
\beaa
E[f\circ V]&=&E[f\circ V\circ U\,\rho(-\delta u)]\\
&=&E[f\,\rho(-\delta u)]\,,
\eeaa
hence $V\mu$ is equivalent to $\mu$ and the corresponding
Radon-Nikodym density is $\rho(-\delta u)$. Let 
$$
D=\{w\in W:\,V\circ U(w)=w\}\,.
$$
Since $D\subset U^{-1}(U(D))$ and by  the hypothesis $\mu(D)=1$ we get
$$
E[1_{U(D)}\circ U]=1\,.
$$
Since $U\mu$ is equivalent to $\mu$ we have also $\mu(U(D))=1$. If
$w\in U(D)$, then $w=U(d)$, for some $d\in D$, hence $U\circ
V(w)=U\circ V\circ U(d)=U(d)=w$, consequently $U\circ V=I_W$
$\mu$-almost surely and $V$ is the two-sided inverse of
$U$. Evidently, together with the absolute continuity of $V\mu$, this
implies that $V$ is of the form $V=I_W+v$, with $v:W\to H$.  Moreover,
$\dot{u}=\dot{v}\circ U$, hence the right hand side is adapted. We can
assume that all these processes are uni-dimensional (otherwise we
proceed component wise). Let
$\dot{v}^n=\max(-n,\min(\dot{v},n))$. Then $\dot{v}^n\circ U$ is
adapted. Let $H\in L^2(dt\times d\mu)$ be an adapted process. Using
the Girsanov theorem:
\beaa
E\left[\rho(-\delta u)\int_0^1\dot{v}^n_s\circ U\,H_s\circ U
ds\right]&=&E\left[\int_0^1\dot{v}^n_sH_sds\right]\\
&=&E\left[\int_0^1E[\dot{v}^n_s|\calF_s]H_sds\right]\\
&=&E\left[\rho(-\delta u)\int_0^1E[\dot{v}^n_s|\calF_s]\circ U\,H_s\circ U
ds\right]\,.
\eeaa
Consequently
$$
E[\dot{v}^n_s|\calF_s]\circ U=\dot{v}^n_s\circ U\,,
$$
almost surely. Since $U\mu$ is equivalent to $\mu$, it follows that 
$$
E[\dot{v}^n_s|\calF_s]=\dot{v}^n_s
$$
almost surely, hence $\dot{v}^n$ and also $\dot{v}$ are adapted. It is
now clear that $(V(t),\,t\in [0,1])$ is a strong solution of
(\ref{sde}). The uniqueness follows from the fact that, any strong
solution of (\ref{sde}) would be a right  inverse to $U$, since $U$ is
invertible, then this solution is equal to $V$.
\qed

\subsection{Carleman inequality}
\noindent
In the sequel we shall use the inequality of T. Carleman which says
that (cf. \cite{Carl} or \cite{D-S}, Corollary XI.6.28)
$$
\|\dett(I_H+A)(I_H+A)^{-1}\|\leq \exp\frac{1}{2}\left(\|A\|_2^2+1\right)\,,
$$
for any Hilbert-Schmidt operator $A$, where the left hand side is the
operator norm, $\dett(I_H+A)$ denotes the modified Carleman-Fredholm
determinant and $\|\cdot\|_2$ denotes the Hilbert-Schmidt norm. Let us
remark that if $A$ is a quasi-nilpotent operator, i.e., if the
spectrum of $A$ consists of zero only, then $\dett(I_H+A)=1$, hence in
this case the Carleman inequality reads 
$$
\|(I_H+A)^{-1}\|\leq
\exp\frac{1}{2}\left(\|A\|_2^2+1\right)\,.
$$
This case happens when $A$ is equal to the Sobolev derivative of some
$u\in \DD_{p,1}(H)$ whose  drift $\dot{u}$ is
adapted to the filtration $(\calF_t,\,t\in  [0,1])$, cf. \cite{ASU,ASU-1}.

\section{A sufficient condition for invertibility}

\label{suff-cond-1}
In the sequel, for a given $u\in \DD_{2,0}(H)$ adapted we shall denote  $e^{-\tau}P_\tau u$ and
$e^{-\ka}P_\ka u$ by $u_\tau$ and $u_\ka$ respectively, the
reason for that is simply the identity $P_\tau\delta u=\delta u_\tau$
is more practical for controlling the Girsanov exponential. Besides we
shall suppose that $u$ satisfies always the condition that
$\rho(-\delta u)$ is a probability density with respect to $\mu$. We then have 
\begin{lemma}
Assume that $u$ is adapted and that
\begin{equation}
\label{tau-con}
E\left[\exp\left(-\delta u-\frac{1}{2}|u|_H^2\right)\right]=E[\rho(-\delta u)]=1\,.
\end{equation}
Then we have 
$$
E\left[\exp\left(-\la\delta
    u_\tau-\frac{\la^2}{2}|u_\tau|_H^2\right)\right]=E[\rho(-\la\delta
u_\tau)]=1\,,
$$
for any $\la,\tau\in [0,1]$.
\end{lemma}
\proof
Define the stopping time 
$$
T_n=\inf\left(t:\int_0^t|\dot{u}_s|^2ds>n\right)\,,
$$
and let $u^n(t)=u(t\wedge T_n)$. Let $u^n_\tau=e^{-\tau}P_\tau u^n$,
then $\lim \rho(- \delta\la  u^n_\tau)=\rho(-\delta u)$ in
probability when $n\to\infty,\,\la\to 1$ and $\tau\to 0$. Besides they
have the constant expectation which is one. Hence $\{\rho(-\delta\la
u^n_\tau):\,n\geq 1,\,\la\in [0,1],\,\tau\in[0,1]\}$ is uniformly
integrable. Consequently its subset $\{\rho(-\delta\la
u^n_\tau):\,n\geq 1\}$ is also uniformly integrable and this completes
the proof.
\qed

\begin{remarkk}
We note also that the  hypothesis \ref{tau-con} is satisfied as
soon as $u$ satisfies either  Novikov or Kazamaki condition,
cf. \cite{I-W}.
\end{remarkk}

\noindent
Theorem \ref{jacobi-thm} implies then that $U_\tau=I_W+u_\tau$ and
$U_\ka=I_W+u_\ka$ are invertible and their inverses are of the form 
$V_\tau=I_W+v_\tau,\,V_\ka=I_W+v_\ka$ respectively. Moreover $v_\tau$
and $v_\ka$ are $H$-valued and adapted.
For $\alpha \in[0,1]$,
let 
$$
u_\alpha^{\tau,\ka}=\alpha u_\tau +(1-\alpha)u_\ka\,.
$$
Then $u_\alpha^{\tau,\ka}$ is again an $H-C^1$-map, it is adapted and
it  inherits all the integrability properties of $u$. 
 Consequently the
map $U_\alpha^{\tau,\ka}$, defined by
$$
w\to w+u_\alpha^{\tau,\ka}(w)
$$
is invertible and its inverse is of the form
$V_\alpha^{\tau,\ka}=I_W+v_\alpha^{\tau,\ka}$ where  $v_\alpha^{\tau,\ka}$
is adapted, $H$-valued, $H-C^1$ and it satisfies the relation 
$$
v_\alpha^{\tau,\ka}=-u_\alpha^{\tau,\ka}\circ V_\alpha^{\tau,\ka}
$$
a.s. Moreover, $v_\tau=v_1^{\tau,\ka}$ and
$v_\ka=v_0^{\tau,\ka}$. 

\noindent
We need the following result:
\begin{lemma}
\label{tech-lemma}
The mapping $\alpha\to v_\alpha^{\tau,\ka}$ is almost surely
continuously differentiable on the interval $(0,1)$.
\end{lemma}
\proof
Define the partial map $t_w^\alpha:H\to H$ as 
$$
t_w^\alpha(h)=h+u_\alpha^{\tau,\ka}(w+h)
$$
for $w\in W$ fixed. Note that from the $H-C^1$-property of $u_\tau$
and $u_\ka$, this map is $C^1$ on $H$ for all $w\in W$ outside a set
of zero capacity. Define the map $\ga$ from $(0,1)\times H$ to itself
as $\ga(\alpha,h)=(\alpha, t_w^\alpha(h))$. Then the differential of
$\ga$ has a Carleman-Fredholm determinant which is equal to
one. Consequently it is invertible as an operator, hence the inverse function theorem
implies the existence of a differentiable inverse $\ga^{-1}$ of
$\ga$. Besides this inverse can be written as
$\ga^{-1}(\alpha,h)=(\alpha,s_w^\alpha(h))$ where $s_w^\alpha$
satisfies the identity 
$$
t_w^\alpha\circ s_w^\alpha=I_H\,, 
$$
where $\alpha\to s_w^\alpha(h)$ is $C^1$ on $(0,1)$. It is easy to see
that $v_\alpha^{\tau,\ka}(w)=s_w^\alpha(0)$ and this completes the
proof.

\qed

\noindent
Hence, due to Lemma \ref{tech-lemma}  we have the following obvious  relation
$$
v_\tau -v_\ka=\int_0^1\frac{dv_\alpha^{\tau,\ka}}{d\alpha}d\alpha\,.
$$

\begin{theorem}
\label{convex-thm}
We have the following inequality:
\begin{equation}
\label{major-eqn}
E[|v_\tau -v_\ka|_H]\leq E\left[|u_\tau-u_\ka|_H\int_0^1
  \exp\frac{1}{2}(\|\nabla u_\alpha^{\tau,\ka}\|_2^2+1)\rho(-\delta
  u^{\tau,\ka}_\alpha) d\alpha\right]
\end{equation}
\end{theorem}
\proof
From Lemma \ref{tech-lemma}, it follows immediately via the chain rule
that 
$$
\frac{dv_\alpha^{\tau,\ka}}{d\alpha}=-(u_\tau-u_\ka)\circ
V_\alpha^{\tau,\ka}-\nabla u_\alpha^{\tau,\ka}\circ
V_\alpha^{\tau,\ka}\,\,\frac{dv_\alpha^{\tau,\ka}}{d\alpha}\,.
$$
Therefore
$$
\frac{dv_\alpha^{\tau,\ka}}{d\alpha}=-\left[(I_H+\nabla
  u_\alpha^{\tau,\ka})^{-1}(u_\tau-u_\ka)\right]\circ
V_\alpha^{\tau,\ka}\,.
$$
Since
$$
\frac{dV_\alpha^{\tau,\ka}\mu}{d\mu}=\rho(-\delta
u_\alpha^{\tau,\ka})\,,
$$
we have
\beaa
E[|v_\tau -v_\ka|_H]&\leq&E\int_0^1\left|\frac{dv_\alpha^{\tau,\ka}}{d\alpha}\right|_Hd\alpha\\
&=&E\int_0^1|(I_H+\nabla
  u_\alpha^{\tau,\ka})^{-1}(u_\tau-u_\ka)|_H\circ
V_\alpha^{\tau,\ka}\,d\alpha\\
&=&E\int_0^1|(I_H+\nabla
  u_\alpha^{\tau,\ka})^{-1}(u_\tau-u_\ka)|_H\,\rho(-\delta
  u^{\tau,\ka}_\alpha) d\alpha\,.
\eeaa
Remarking that  $\nabla u_\alpha^{\tau,\ka}$ is quasi-nilpotent and
applying the Carleman inequality in the last line of the above inequalities, we get 
$$
E[|v_\tau -v_\ka|]\leq E\left[|u_\tau-u_\ka|_H\int_0^1
  \exp\frac{1}{2}(\|\nabla u_\alpha^{\tau,\ka}\|_2^2+1)\rho(-\delta
  u^{\tau,\ka}_\alpha) d\alpha\right]
$$
and this completes the proof.
\qed

\begin{theorem}
\label{nice-thm}
Assume that $u\in \DD_{p,1}(H)$ for some $p>1$ and that it  is adapted with $E[\rho(-\delta
u)]=1$. Suppose moreover  that $u$ satisfies  the following condition:
$$
E\left[\exp q\left(\frac{1}{2}\|\nabla u\|_2^2-\delta u\right)\right]<\infty\,,
$$
where $p^{-1}+q^{-1}=1$. Then $U=I_W+u$ is almost surely invertible.
\end{theorem}
\proof
From Theorem \ref{convex-thm}, using the H\"older inequality  we have 
\beaa
\lefteqn{E[|v_\tau -v_\ka|_H]\leq E\left[|u_\tau-u_\ka|_H\int_0^1
  \exp\frac{1}{2}(\|\nabla u_\alpha^{\tau,\ka}\|_2^2+1)\rho(-\delta
  u^{\tau,\ka}_\alpha) d\alpha\right]}\\
&&\lefteqn{\leq E\left[|u_\tau-u_\ka|_H\int_0^1\exp\left(\frac{\alpha}{2}\|\nabla
u_\tau\|_2^2-\alpha\delta u_\tau+\frac{1-\alpha}{2}\|\nabla
u_\ka\|_2^2-(1-\alpha)\delta u_\ka\right)d\alpha\right]}\\
&\leq& E[|u_\tau-u_\ka|_H^p]^{1/p}\\
&&\,\times\,\left[E\int_0^1\exp q\left(\frac{\alpha}{2}\|\nabla
u_\tau\|_2^2-\alpha\delta u_\tau+\frac{1-\alpha}{2}\|\nabla
u_\ka\|_2^2-(1-\alpha)\delta u_\ka\right)d\alpha\right]^{1/q}\\
&\leq&E[|u_\tau-u_\ka|_H^p]^{1/p}\\
&&\times\,\left(\int_0^1E\left[\exp q\left(\frac{1}{2}\|\nabla
u_\tau\|_2^2-\delta u_\tau\right)\right]^\alpha\,\times\,E\left[\exp q\left(\frac{1}{2}\|\nabla
u_\ka\|_2^2-\delta u_\ka\right)\right]^{1-\alpha}d\alpha\right)^{1/q}
\eeaa
From the Jensen inequality and from the relation
$$
 \delta u_\tau=P_\tau \delta u\,,
$$ 
we obtain
$$
E\left[\exp q\left(\frac{1}{2}\|\nabla u_\tau\|_2^2-\delta
    u_\tau\right)\right]\leq E\left[\exp q\left(\frac{1}{2}\|\nabla
u\|_2^2-\delta u \right)\right]\,.
$$
Consequently 
$$
E[|v_\tau -v_\ka|_H]\leq E[|u_\tau-u_\ka|_H^p]^{1/p}E\left[\exp q\left(\frac{1}{2}\|\nabla
u\|_2^2-\delta u \right)\right]^{1/q}\rightarrow 0
$$
 since $E[|u_\tau-u_\ka|_H^p]\to 0$ as $\ka,\tau\to 0$ and this
 implies the existence of some adapted  $v:W\to H$ which is the limit in
 $L^1(\mu,H)$ of $(v_\tau,\,\tau\in (0,1))$. To complete the proof we
 have to show that $v\circ U=-u$ and $u\circ V=-v$ almost surely,
 where $V=I_W+v$. For $c>0$, we have
\begin{eqnarray}
\label{eqn-11}
\mu\left\{|v_\tau\circ U_\tau-v\circ U|_H>c\right\}&\leq&\mu\left\{|v_\tau\circ
U_\tau-v\circ U_\tau|_H>\frac{c}{2}\right\}\nonumber\\
&&+\mu\left\{|v\circ U_\tau-v\circ U|_H>\frac{c}{2}\right\}\nonumber\\
&=&E\left[\rho(-\delta v_\tau)1_{\{|v_\tau-v|_H>c/2\}}\right]\\
&&+\mu\left\{|v\circ U_\tau-v\circ U|_H>\frac{c}{2}\right\}\label{eqn-12}\,
\end{eqnarray}
Since 
$$
E[\rho(-\delta v_\tau)\log \rho(-\delta
v_\tau)]=\frac{1}{2}E[|u_\tau|_H^2]\,,
$$
the set $(\rho(-\delta v_\tau),\,\tau\in [0,1])$ is uniformly
integrable, hence  the first term  (\ref{eqn-11}) can be made arbitrarily
small by the convergence of $v_\tau\to v$ in probability. Moreover, we
know that $(\rho(-\delta u_\tau),\,\tau\in [0,1])$ converges in
probability to $\rho(-\delta u)$ and they have all the same
expectation which is equal to one. Consequently the set
$(\rho(-\delta u_\tau),\,\tau\in [0,1])$ is also uniformly
integrable. To control the term (\ref{eqn-12}), recall that, by the
Lusin theorem, given any $\eps>0$, there exists a compact set $K_\eps$
in $W$ such that $\mu(K_\eps)>1-\eps$ and that the restriction of $v$
to $K_\eps$ is uniformly continuous. Therefore
\begin{eqnarray}
\lefteqn{\mu\left\{|v\circ U_\tau-v\circ U|_H>\frac{c}{2}\right\}}\nonumber\\
&\leq&
\mu\left\{|v\circ U_\tau-v\circ U|_H>\frac{c}{2},\,U_\tau\in
  K_\eps,\,U\in K_\eps\right\}\label{eqn-13}\\
&&+\mu\{U_\tau\in K_\eps^c\}+\mu\{U\in K_\eps^c\} \label{eqn-14}
\end{eqnarray}
The last two terms (\ref{eqn-14})  can be made arbitrarily small (uniformly
w.r. to $\tau$) by the uniform integrability of $(\rho(-\delta
u_\tau),\,\tau\in [0,1])$. To control the term (\ref{eqn-13}), let
$\beta>0$ be arbitrary. Then 
\begin{eqnarray}
\lefteqn{\mu\left\{|v\circ U_\tau-v\circ U|_H>\frac{c}{2},\,U_\tau\in
  K_\eps,\,U\in K_\eps\right\}}\nonumber\\
&\leq&\mu\left\{|v\circ U_\tau-v\circ U|_H>\frac{c}{2},\,U_\tau\in
  K_\eps,\,U\in K_\eps,\,\|U_\tau-U\|>\beta\right\}\label{eqn-15}\\
&&+\mu\left\{|v\circ U_\tau-v\circ U|_H>\frac{c}{2},\,U_\tau\in
  K_\eps,\,U\in K_\eps,\,\|U_\tau-U\|\leq\beta\right\}\label{eqn-16}
\end{eqnarray}
where $\|\cdot\|$ denotes the norm of $W$. Since $v$ is uniformly
continuous on $K_\eps$, the term (\ref{eqn-16}) can be made
arbitrarily small by choosing $\beta$ small enough and the term
(\ref{eqn-15}) is bounded by 
$$
\mu\{ \|U_\tau-U\|>\beta\}
$$
which can be made arbitrarily small by choosing $\tau$ small enough
and this proves the relation $v\circ U=-u$ which implies that $V\circ
U=I_W$ almost surely. To prove $u\circ V=-v$, recall that 
$$
\frac{dV_\tau\mu}{d\mu}=\rho(-\delta u_\tau)
$$
and as we have indicated above $(\rho(-\delta u_\tau),\,\tau\in
[0,1])$ is uniformly integrable. Hence we can repeat the same
reasoning as  above
by interchanging $u$ and $v$ in the above lines and this completes the
proof.
\qed

\begin{corollary}
\label{holder-cor}
Assume that $u\in \DD_{p,1}(H)$ is adapted. If $u$ satisfies the
following condition
$$
E\left[\exp\left(q\|\nabla u\|_2^2+2q^2|u|_H^2\right)\right]<\infty\,,
$$
then $U=I_W+u$ is almost surely invertible.
\end{corollary}
\proof
Let $\eps>1$, we have, using the H\"older inequality  
\begin{align*}
E\left[\exp q\left(\frac{1}{2}\|\nabla u\|_2^2-\delta u\right)\right]
    &=E\left[\exp
\left(q\frac{1}{2}\|\nabla u\|_2^2-q\delta
u-q^2\frac{1+\eps}{2\eps}|u|_H^2+q^2\frac{1+\eps}{2\eps}|u|_H^2\right)\right]\\
&\leq E\left[\exp\left(\frac{(1+\eps)q}{2}\|\nabla
u\|_2^2+q^2\frac{(1+\eps)^2}{2\eps}|u|_H^2\right)\right]^{1/1+\eps}\\
&\qquad\qquad \times \, E\left[\exp\left(-\frac{(1+\eps)}{\eps}q\delta
u-q^2\frac{(1+\eps)^2}{2\eps^2}|u|_H^2\right)\right]^{\eps/1+\eps}\\
&\leq E\left[\exp\left(\frac{(1+\eps)q}{2}\|\nabla
u\|_2^2+q^2\frac{(1+\eps)^2}{2\eps}|u|_H^2\right)\right]^{1/1+\eps}\,,
\end{align*}
since the expectation of  third line is upperbounded by one. The proof
follows when we take $\eps=1$ for which the last line attains its
minimum with respect to $\eps>0$.
\qed

\label{log-sobolev}
\begin{theorem}
\label{log-sob}
Assume that $L$ is a probability density with respect to $\mu$ which
has an It\^o  representation 
$$
L=\rho(-\delta u)=\exp\left[-\delta u-\frac{1}{2}|u|_H^2\right]\,,
$$
where $u\in \DD_{2,1}(H)$ is adapted. If $\|\nabla u\|_2\in
L^\infty(\mu)$, then the measure $\nu$, defined as 
$$
d\nu=L\,d\mu
$$
satisfies the logarithmic Sobolev inequality, i.e., 
$$
E_\nu\left[f^2\log\frac{f^2}{E_\nu[f^2]}\right]\leq K E_\nu[|\nabla f|_H^2]\,,
$$
for any cylindrical Wiener function $f$, where 
$$
K=2\left\|\exp\left(1+\|\nabla
    u\|_2^2\right)\right\|_{L^\infty(\mu)}\,.
$$
\end{theorem}
\proof 
We shall use a reasoning analogous to that of \cite{ASU-2}. Let
  $\alpha$ be a positive, smooth function of
compact support with $\alpha(0)=1$ and $|\alpha'(t)|\leq c$ for any
$t\geq 0$. Define $\dot{u}^n_t$ and $u^n$ as
$$
\dot{u}_t^n=\alpha(\frac{1}{n}\dot{u}_t)\,\dot{u}_t,\,u^n(t)=\int_0^t\dot{u}^n_s
ds\,.
$$
$u^n$ is bounded and 
$$
\|\nabla u^n\|_2^2\leq 2(1+c^2)\|\nabla u\|_2^2\in L^\infty(\mu)\,.
$$
From Theorem \ref{nice-thm},  $U^n=I_W+u^n$ is a.s. invertible and its
inverse $V^n$ is of the form $I_W+v^n$ such that  $v^n$ is in 
$\DD_{2,1}(H)$ and adapted. Consequently 
$$
\frac{dV^n\mu}{d\mu}=\rho(-\delta u^n)\,.
$$
Let $\nu_n$ be the probability measure defined as $d\nu_n=\rho(-\delta
u^n)d\mu$. Using the log-Sobolev inequality of L. Gross for
$\mu$, cf. \cite{Gross},  and the Carleman inequality  we get
\beaa
E_{\nu_n}\left[f^2\log\frac{f^2}{E_{\nu_n}[f^2]}\right]&=&E\left[(f\circ
  V^n)^2\log\frac{f^2\circ V^n}{E_{\nu_n}[f^2]}\right]\\
&\leq&2\,E\left[|\nabla (f\circ V^n)|_H^2\right]\\
&\leq&2\,E\left[|\nabla f\circ V^n|_H^2\|I_H+\nabla v^n\|^2\right]\\
&=&2\,E\left[\rho(-\delta u^n)|\nabla f|_H^2\|(I_H+\nabla
  u^n)^{-1}\|^2\right]\\
&\leq&\,E\left[\rho(-\delta u^n)|\nabla f|_H^2\exp\left(1+\|\nabla
    u^n\|_2^2\right)\right]\\
&\leq&\,2\left\|\exp\left(1+\|\nabla
  u^n\|_2^2\right)\right\|_{L^\infty(\mu)}E\left[\rho(-\delta u^n)|\nabla
f|_H^2\right]\\
&\leq&2E_{\nu_n}[|\nabla f|_H^2]\exp\left(1+2(1+c^2)\|\|\nabla u\|_2\|^2_{L^\infty(\mu)}\right)\,.
\eeaa
To complete the proof,  take first  the limit of this inequality
as $n\to \infty$ and remark that $(\rho(-\delta u^n),\,n\in
\N)$ is uniformly integrable. Finally it suffices  to take  the
infimum of the right hand side with respect to $c>0$.
\qed

\begin{theorem}
\label{lip-thm}
Under the hypothesis of Theorem \ref{log-sob}, the mapping $U=I_W+u$
is almost surely invertible.
\end{theorem}
\proof
With the notations of Theorem \ref{nice-thm}, we have 
$$
E[|v_\tau-v_\kappa|_H]\leq E\left[|u_\tau-u_\kappa|_H\int_0^1\rho(-\delta u_\alpha^{\tau,\kappa})d\alpha\,\right]\,.
$$
Let us denote $\rho(-\delta u^{\tau,\kappa}_\alpha)$ by
$\rho^{\tau,\kappa}_\alpha$ and let $\nu_\alpha^{\tau,\kappa}$ be the
measure whose Radon-Nikodym derivative with respect to $\mu$ is given
by $\rho^{\tau,\kappa}_\alpha$. We have
\beaa
E[\rho^{\tau,\kappa}_\alpha
|u_\tau|_H^2]&\leq&2E[\rho^{\tau,\kappa}_\alpha
|u_\tau-E_{\nu_\alpha^{\tau,\kappa}}[u_\tau]|_H^2]+
2|E_{\nu_\alpha^{\tau,\kappa}}[u_\tau]|_H^2\\
&=&I_{\tau,\kappa,\alpha}+II_{\tau,\kappa,\alpha}\,.
\eeaa
From Theorem \ref{log-sob} and from the fact that logarithmic Sobolev
inequality implies the Poincar\'e inequality, the first terms at the
right hand side of the above inequality is bounded:
\begin{equation}
\label{bound-1}
I_{\tau,\kappa,\alpha}\leq 2 C\,E_{\nu_\alpha^{\tau,\kappa}}[\|\nabla u_\tau\|_2^2]\leq
2C\|\nabla u\|_{L^\infty(\mu,H\otimes H)}^2\,,
\end{equation}
where $C$ is independent of $\alpha,\kappa$ and $\tau$; in fact it is
the constant of the logarithmic Sobolev inequality of Theorem
\ref{log-sob}. Moreover
$$
u_\tau=(u_\tau-E_{\nu_\alpha^{\tau,\kappa}}[u_\tau])+E_{\nu_\alpha^{\tau,\kappa}}[u_\tau]\,.
$$
From the inequality (\ref{bound-1}), it follows that 
$$
\sup_{\tau,\kappa,\alpha}\nu^{\tau,\kappa}_\alpha(|u_\tau-E_{\nu_\alpha^{\tau,\kappa}}[u_\tau]|_H>c)\to
0
$$
as $c\to\infty$. Besides, from the uniform integrability of
$(\rho_\alpha^{\tau,\kappa}:\,\tau,\kappa,\alpha\in [0,1])$, we also
have 
$$
\sup_{\tau,\kappa,\alpha}\nu^{\tau,\kappa}_\alpha(|u_\tau|_H>c)\leq\sup_{\tau,\kappa,\alpha}E[\rho^{\tau,\kappa}_\alpha
1_{\{|u_\tau|_H>0\}}]\to 0
$$
as $c\to\infty$. Consequently, there exists a constant $c>0$ such that
$$
\sup_{\tau,\kappa,\alpha\in[0,1]}|E_{\nu_\alpha^{\tau,\kappa}}[u_\tau]|_H
\leq c\,.
$$
This implies the uniforme integrability of the family 
$$
\left(|u_\tau-u_\kappa|_H\int_0^1\rho(-\delta
u_\alpha^{\tau,\kappa})d\alpha:\,\tau,\kappa\in[0,1]\right)\,.
$$
Since it converges converges already in probability to zero, the
convergence in $L^1(\mu)$ holds also. The rest of the proof follows
the same lines as the proof of Theorem \ref{nice-thm}, hence the it is
completed.
\qed

\begin{remarkk}
\label{remark-1}
In terms of the  stochastic differential equations, the question of
finding an inverse to $U=I_W+u$, where $u\in \DD_{2,0}(H)$ is adapted,
amounts to solving the following stochastic differential equation 
\begin{eqnarray}
\label{SDE}
dV_t(w)&=&-\dot{u}_t(V(w))dt+ dW_t(w)\\
V_0(w)&=&0\nonumber\,,
\end{eqnarray}
and this problem  is solved only  under a  Lipschitz hypothesis imposed
to $\dot{u}$, which can be expressed  as follows
\begin{equation}
\label{lip-cond}
\sup_{s\leq t}|\dot{u}_s(w+h)-\dot{u}_s(w+k)|\leq K\,\sup_{s\leq
  t}|h(s)-k(s)|\,,
\end{equation}
$\mu$-a.s., for any $h,k\in H$, $t\in [0,1]$,  where $K$ is a constant. Since 
$$
\sup_{s\leq t}|h(s)-k(s)|\leq |h-k|_H\,,
$$
the Lipschitz condition (\ref{lip-cond}) implies that 
$$
|\nabla\dot{u}_s|_H\leq K\,,
$$
hence 
$$
\|\nabla u\|_2^2=\int_0^1|\nabla \dot{u}_s|_H^2ds \leq K^2
$$
$\mu$-almost surely. 
Therefore, the Lipschitz condition (\ref{lip-cond}) is stronger than
the hypothesis of Theorem \ref{lip-thm}. For example, assume that
the drift $\dot{u}$ has a Sobolev derivative which satisfies 
$$
|\nabla\dot{u}_s|_H\leq K\,s^{-\alpha}\,,
$$
almost surely, where $0\leq \alpha<1/2$. Then the Lipschitz property may
fail although the stochastic differential equation (\ref{SDE}) has a
unique solution by the theorem  since $\|\nabla u\|_2\in L^\infty(\mu)$.
\end{remarkk}

\section{Extensions}
\label{extensions}
In this section we give some variations and extensions  of the results
proven in the last section. We start with 
\begin{theorem}
\label{extension-thm}
Let $u:W\to H$ be adapted with $u\in \DD_{2,0}(H)$ such that
$E[\rho(-\delta u)]=1$. Assume that $(\Om_n,n\geq 1)$ is a measurable covering of $W$ and that
$u=u_n$ a.s. on $\Om_n$ where $u_n:W\to H$ is in $\DD_{2,0}(H)$,
adapted, $E[\rho(-\delta u_n)]=1$  and $U_n=I_W+u_n$ is almost surely invertible with the inverse
denoted by $V_n=I_W+v_n$. Then $U=I_W+u$ is almost surely invertible and the
Radon-Nikodym derivative of $U\mu$ with respect to $\mu$ belongs to
the space $L\log L(\mu)$. 
\end{theorem}
\proof
Without loss of generality we can assume that the sets $(\Om_n, n\geq
1)$ are disjoint. Note also that  by the hypothesis, $U\mu$ is equivalent to $\mu$. 
 We  have  
$$
\frac{dU_n\mu}{d\mu}=\rho(-\delta v_n)\,.
$$
 Then, for any $f\in C_b(W)$ 
\beaa
E[f\circ U\,\rho(-\delta u)]&=&\sum_{n=1}^\infty E[f\circ U_n\,1_{\Om_n}\rho(-\delta u_n)]\\
&=&\sum_{n=1}^\infty E[f\,\,1_{U_n(\Om_n)}]\,.
\eeaa
By the Girsanov theorem we also have 
$$
E[f\circ U\,\rho(-\delta u)]=E[f]\,,
$$
hence 
$$
\sum_{n=1}^\infty 1_{U_n(\Om_n)}=1
$$
almost surely. This means that $(U_n(\Om_n),n\geq 1)$ is an almost
sure partition of $W$. Define $v$ on $U_n(\Om_n)$ as to be $v_n$ and
let $V=I_W+v$. Then $V$ is defined almost everywhere, moreover
\beaa
E[f\circ V]&=&\sum_nE[f\circ V_n\,1_{U_n(\Om_n)}]\\
&=&E[f\,\rho(-\delta u)]\,,
\eeaa
therefore $V\mu$ is equivalent to $\mu$ and $V$ is
well-defined. Evidently, for almost all  $w\in \Om_n$,
$$
V\circ U(w)=V\circ U_n(w)=V_n\circ U_n(w)=w\,,
$$
hence $V\circ U=I_W$ almost surely, i.e., $V$ is a left inverse of
$U$. Since
$$
\mu\left(\bigcup_{n=1}^\infty U_n(\Om_n)\right)=1\,,
$$
and since 
$$
\{w\in W:\,U\circ V(w)=w\}\supset \bigcup_{n=1}^\infty U_n(\Om_n)\,,
$$
we also have
$$
\mu\left(\{w\in W:\,U\circ V(w)=w\}\right)=1
$$
and this completes the proof of the invertibility of $U$. Clearly, the
map $v$ is adapted to the filtration of the Wiener space, hence the
stochastic integral of the drift $(\dot{v}_t,t\in [0,1])$,  with
respect to the Wiener process is well-defined
(using the localization techniques with the help of the stopping
times) and we shall denote its value at $t=1$ by $\delta^0v$. Since
$U$ is adapted, we have 
\beaa
(\delta^0v)\circ U&=&\delta^0(v\circ U)+(v\circ U,u)_H\\
&=&-\delta^0u-|u|_H^2\\
&=&-\delta u-|u|_H^2\,,
\eeaa
hence 
$$
\rho(-\delta^0v)\circ U\,\rho(-\delta u)=1\,,
$$
 and similarly
$$
\rho(-\delta u)\circ V\,\rho(-\delta^0 v)=1
$$
almost surely. Therefore, on the one hand 
$$
\frac{dU\mu}{d\mu}=\rho(-\delta^0v)
$$
and on the other hand
$$
E\left[\rho(-\delta^0v)\log \rho(-\delta^0 v)\right]=\frac{1}{2}E[|u|_H^2]<\infty
$$
by the  hypothesis and the proof is completed.

\qed

\noindent

Here is another situation which is encountered in the applications:
\begin{theorem}
\label{stopping-thm}
Assume that $u\in \DD_{2,0}^a(H)$ with $E[\rho(-\delta u)]=1$. Suppose
that there exists a sequence of stopping times $(T_n,n\geq 1)$
increasing to infinity and a sequence $(u_n,n\geq 1)\subset
\DD_{2,0}^a(H)$ such that 
$$
\dot{u}(t)=\dot{u}_n(t)\,\,\,{\mbox{for }}\,t<T_n(w)\,.
$$
Suppose further that $U_n=I_W+u_n$ are almost surely invertible with
inverse $V_n=I_W+v_n$ for any $n\geq 1$. Then $U=I_W+u$ is almost surely invertible with
inverse $V=I_W+v$, $v\in L^0(\mu,H)$, with $\dot{v}$ adapted and 
$$
\frac{dU\mu}{d\mu}=\rho(-\delta^0v)\,.
$$
\end{theorem}
\proof
For any $n\geq 1$, $(V_n(t),t\in [0,1])$ is the unique solution of the
equation 
$$
V_n(t)=W_t-\int_0^t\dot{u}_n(s,V_n)ds\,.
$$
By hypothesis, for $m\leq n$, we have 
$$
\dot{u}_m(t\wedge T_m,w)=\dot{u}_n(t\wedge T_m,w)=\dot{u}(t\wedge T_m,w)
$$
$dt\times d\mu$-almost surely. Hence 
$$
V_n(t)=W_t-\int_0^t\dot{u}_m(s,V_n)ds
$$
for $t\leq T_m(V_n(w))$. Hence, by the uniqueness
\begin{equation}
\label{local}
V_m(t)=V_n(t)\,\,{\mbox{ for any }}\,t\leq T_m(V_n(w))\,.
\end{equation}
Similarly, for $t\leq T_n(V_m)$ we have 
$$
V_m(t)=W_t-\int_0^t\dot{u}_n(s,V_m)ds\,,
$$
hence again by the uniqueness
\begin{equation}
\label{local-1}
V_m(t)=V_n(t)\,\,{\mbox{ for any }}\,t\leq T_n(V_m(w))\,.
\end{equation}
Consequently, combining the relations (\ref{local}) and
(\ref{local-1}), we get 
\begin{equation}
\label{local-2}
V_m(t)=V_n(t)\,\,{\mbox{ for any }}\,t\leq T_n(V_m(w))\vee
T_m(V_n(w))\,.
\end{equation}
In particular, for any $m$ and for any $n\geq m$, we have  
\begin{equation}
\label{local-3}
V_m(t)=V_n(t)\,\,{\mbox{ for any }}\,t\leq T_m(V_m(w))\,.
\end{equation}
Let $\tilde{T}_n=\sup_{k\leq n}T_k(V_k)$, we have
\beaa
\mu\{T_n(V_n)>t\}&=&E[1_{\{T_n>t\}}\rho(-\delta u_n)]\\
&=&E[1_{\{T_n>t\}}E[\rho(-\delta u_n)|\calF_{T_n}]]\\
&=&E[1_{\{T_n>t\}}\rho(-\delta u)]\to 1\,,
\eeaa
therefore $(\tilde{T}_n,n\geq 1)$ increases to infinity. Define now
$(V(t),t\in [0,1])$ as follows: 
$$
V(t)=V_n(t) {\mbox{ if }}\,t\leq \tilde{T}_n(w)\,.
$$
In fact, for any $t\in [0,1]$ there exists some $k\leq n$ such that $t\leq
T_k(V_k)$, hence by the relation (\ref{local-3}), 
$$
V_k(t)=V_n(t)=V_{n+l}(t)\,\,{\mbox{ for }}t\leq T_k(V_k)\,,
$$
for any $l\geq 1$. In particular, for $m\leq n$ and $t\leq
\tilde{T}_m(w)$,  we have 
$$
V_m(t)=V_n(t)\,,
$$
hence $V$ is well-defined.
Moreover, $V$ can be written as 
$$
V(t)=W_t+\int_0^t\dot{v}_sds\,,
$$
with $\dot{v}\in L^0(\mu,L^2([0,1]))$ adapted. Besides, for $t\leq
{\tilde{T}}_n$, since $\dot{u}_n$ is adapted
\beaa
V(t)&=&W_t-\int_0^t\dot{u}_n(s,V_n)ds\\
&=&W_t-\int_0^t\dot{u}_n(s,V)ds\,,
\eeaa
moreover, from the hypothesis
$$
\dot{u}_n(s,V)=\dot{u}(s,V)\,{\mbox{ for }}s<T_n(V)\,.
$$
Since $V\mu$ is absolutely continuous with respect to $\mu$,
$(T_n(V),n\geq 1)$ incerases to infinity almost surely. Consequently,
for any $t\in [0,1]$, 
$$
V(t)=W_t-\int_0^t\dot{u}(s,V)ds\,,
$$
almost surely. This means that $U\circ V=I_W$ almost surely, in other
words $V$ is a right inverse to $U$.

Let us show now that the mapping  $V$ constructed above is also a left
inverse: we have 
$$
v_m(t,w)=v_n(t,w)=v(t,w)\,{\mbox{ if }} t\leq T_m(V_m),\,m\leq n\,.
$$
Hence, for $t\leq T_m(w)$ and  by the adaptedness of $\dot{v}$,
$$
v_m(t,U_m(w))=v_n(t,U_m(w))=v(t,U_m(w))=v(t,U(w))\,.
$$
By the hypothesis, for $t\leq T_m(w)$, we also have
$$
v_m(t,U_m(w))=-u_m(t,w)=-u(t,w)\,.
$$ 
Since $(T_m,m\geq 1)$ increases to infinity, we obtain
$$
v\circ U+u=0
$$ almost surely. This implies that $V$ is also a left inverse and
that 
$$
\rho(-\delta^0v)\circ U\,\,\rho(-\delta u)=1
$$
almost surely, hence $E[\rho(-\delta^0v)]=1$ and $\rho(-\delta^0v)$ is
the Radon-Nikodym derivative of $U\mu$ with respect to $\mu$.
\qed

\begin{corollary}
\label{bdd-k-cor}
Assume that $u\in \DD_{2,k}(H)$ is adapted  with $k\geq 1$ such that
$E[\rho(-\delta u)]=1$ and that  
$$
\|\nabla^k u\|_{H^{\otimes(k+1)}}\in L^\infty(\mu)\,.
$$
Assume moreover that $H^{\otimes i}$-valued process $t\to
\nabla^i\dot{u}_t$ is almost surely continuous for $0\leq i\leq k-1$.  
Then the mapping $U=I_W+u$ is almost surely invertible.
\end{corollary}
\proof
Let $\theta_n$ be a smooth  function on $\R$ which is equal to one on
$[0,n]$ and zero on the complement of $[-1,n+1]$.
Define $u_n$, $n\geq 1$ as 
$$
\dot{u}_n(t)=\dot{u}(t)\,\prod_{i=0}^{k-1}\theta_n\left(\|\nabla^i
 \dot{u}_t\|^2_{H^{\otimes i}}\right)
$$
for $t\in [0,1]$. Then it is easy to see that $\|\nabla u_n\|_2\in
L^\infty(\mu)$, hence, from Theorem  \ref{lip-thm},  $U_n=I_W+u_n$ is
almost surely invertible for any $n\geq 1$. Define the stopping times
$(T_n,n\geq 1)$ as 
$$
T_n=\inf\left(\sum_{i=0}^{k-1}\|\nabla^i\dot{u}_t\|_{H^{\otimes
    i}}>n\right)
$$
By the continuity hypothesis, $(T_n,n\geq 1)$ increases to infinity,
besides for $t<T_n(w)$, $\dot{u}(t,w)=\dot{u}_n(t,w)$, hence the proof
follows from Theorem \ref{stopping-thm}.
\qed

\noindent
In Theorem \ref{extension-thm}, we have supposed that $u=u_n$ on a set
$\Om_n$, where $u_n$ is also adapted. However, we can construct easily
examples where $u_n$ is not adapted but still $I_W+u_n$ is 
invertible and  equal to $I_W+u$ almost surely on $\Om_n$ such that the
union of the sets $(\Om_n,\,n\geq 1)$ is equal to $W$ almost
surely. In such a situation the hypothesis of Theorem
\ref{extension-thm} are not satisfied.
 To study this kind of situations, we need to define some
more regularity concepts which are studied in detail in \cite{BOOK}:
\begin{definition}
\label{hc-1-defn}
 Let $X$ be a separable Hilbert space, then 
\begin{enumerate}
\item  a measurable map $\xi:W\to X$ is
called $H-C^1_{\rm loc}$ if there exists a measurable $q:W\to \R_+$,
$q>0$ a.s., such that the map $h\to \xi(w+h)$ is a $C^1$-map on the
ball $\{h\in H:\,|h|_H<q(w)\}$.
\item $\xi$ is called representable by locally $H-C^1$-functions or
  $RH-C^1_{\rm loc}$ in short, if there is a sequence of measurable sets
  $(B_n,n\geq 1)$ whose union is of full measure and a sequence of
  $H-C^1_{\rm loc}$-functions $(u_n,n\geq 1)$ such that $u=u_n$ on $B_n$
  almost surely.
\end{enumerate}
\end{definition}
Let us recall  Theorem 3.5.3 of \cite{BOOK} which is valid for not
neccessarily adapted perturbations of identity:
\begin{theorem}
\label{book-thm}
If $u:W\to H$ is $RH-C^1_{\rm loc}$, then, there exists a set $\tilde{W}$ of
full-measure such that, for any $f,\,g\in C_b(W)$, one has 
$$
E[f\circ U\,|\La_u|\,g]=E\left[f(w)\sum_{y\in
    U^{-1}\{w\}}g(y)\,1_{M\cap \tilde{W}}(y)\right]\,,
$$
where $M=\{w\in W:\,\dett(I_H+\nabla u(w))\neq 0\}$ and 
$$
\La_u=\dett(I_H+\nabla u)\exp\left(-\delta
  u-\frac{1}{2}|u|_H^2\right)\,.
$$
In particular,
the multiplicity of $U$ on the set $M\cap\tilde{W}$ is almost surely
(atmost) countable.
\end{theorem}

The next theorem answers to the question that we have asked above:
\begin{theorem}
\label{anticip-thm}
Let $u\in \DD_{2,1}(H)$ be adapted and assume $E[\rho(-\delta
u)]=1$. Suppose that there exists $(\Om_n,n\geq 1)\subset \calB(W)$,
whose union is of full measure and a sequence $(u_n,n\geq 1)$ of
$RH-C^1_{\rm loc}$-functions such that $u=u_n$ almost surely on
$\Om_n$ for any $n\geq 1$. Then $U=I_W+u$ is almost surely invertible.
\end{theorem}
\proof 
Without loss of generality, we can assume that the sets $(\Om_n, n\geq
1)$ are disjoint. Then the following change of variables formula is a
consequence of  Theorem \ref{book-thm} and of  the
fact that $u$ is adapted:  
$$
E[f\circ U\,\rho(-\delta u)\,g]=\sum_{n=1}^\infty E\left[f(w)\sum_{y\in
    U_n^{-1}\{w\}}g(y)\,1_{\Om_n\cap \tilde{W}_n}(y)\right]\,.
$$
In particular, taking $g=1$, we see that 
$$
\sum_{n=1}^\infty \sum_{y\in U_n^{-1}\{w\}} 1_{\Om_n\cap{W}_n}(y)=1
$$
almost surely. Let us denote the set $\Om_n\cap{W}_n$ by $\Om_n'$ and
the double sum above by $N_n(w,\Om_n')$. Since each $N_n$ is an
integer and since their sum is equal to one almost surely, we should
have $N_n(w,\Om_n')\in \{0,1\}$ almost surely. Let 
$$
\tilde{\Om}_n=\{w:\,N_n(w,\Om_n')=1\}\,.
$$
If $w\in \tilde{\Om}_n$, then $N_n(w,\Om_n')=1$, i.e., the cardinal of
the set, denoted by $|U_n^{-1}\{w\}\cap \Om_n'|$ is equal to
one. Consequently, there exists a unique $y\in \Om_n'$ such that
$U_n(y)=w$. This means that $U_n:\Om_n'\to \tilde{\Om}_n$ is
surjectif. Denote the map $w\to y$ by $V_n(w)$, hence
$V_n(\tilde{\Om}_n)\subset\Om_n'$. Define $V$ on
$\cup_n\tilde{\Om}_n$ as $V=V_n$ on $\tilde{\Om}_n$. Since the sets
$\Om_n'$ and $\tilde{\Om}_n$ are measurable, $V$ is measurable with
respect to the completed Borel  sigma algebra of $W$. Taking $g=1$ in
the change of variables formula, we get 
\beaa
E[g\,\rho(-\delta u)]&=&E\left[\sum_n\sum_{y\in
    U_n^{-1}\{w\}}g(y)1_{\Om_n'}(y)\right]\\
&=&E\left[\sum_n 1_{\tilde{\Om}_n}(w)g\circ V_n(w)\right]\\
&=&E[g\circ V]\,.
\eeaa
This implies in particular that the measure $V(\mu)$ is equivalent to
$\mu$. To show that $V$ is also a left inverse, choose any two
$f,\,g\in C_b(W)$. Then 
\beaa
E\left[f\circ V\circ U\,\rho(-\delta u)\,g\right]&=&\sum_nE\left[f\circ V\sum_{y\in U_n^{-1}\{w\}}g(y)1_{\Om_n'}(y)\right]\\
&=&\sum_nE\left[f\circ V(w)\,1_{\tilde{\Om}_n}(w)g\circ V_n(w)\,1_{\Om_n'}\circ V_n(w)\right]\\
&=&\sum_nE\left[f\circ V(w)\,1_{\tilde{\Om}_n}(w)g\circ V_n(w)\right]\\
&=&\sum_nE\left[f\circ V(w)\,1_{\tilde{\Om}_n}(w)g\circ V(w)\right]\\
&=&E[f\circ V\,g\circ V]\\
&=&E[f\,g\,\rho(-\delta u)]\,,
\eeaa
where the second line follows from the fact that the sum on the set
$U_n^{-1}\{w\}$ is zero unless $w\in \tilde{\Om}_n$, in which case
$1_{\Om_n'}\circ V_n(w)=1$ since $V_n(\tilde{\Om}_n)\subset\Om_n'$  by
the construction of $V_n$. 
Consequently $V\circ U=I_W$ $\mu$-almost surely, hence $V$ is a two
sided inverse, it is of the form $V=I_W+v$ and  $v:W\to H$  is adapted to the (completed) filtration of the
Wiener space. Moreover $\delta v$ is well-defined as local martingales
final value (using the stopping techniques), in particular the Radon-Nikodym density of $U\mu$ with
respect to $\mu$ is $\rho(-\delta v)$. 
\qed

 
{\footnotesize} 
\begin{tabular}{ll} 
A. S. \"Ust\"unel  & M. Zakai,\\ 
ENST, Paris& Technion, Haifa\\
Dept. Infres, 46, rue Barrault & Dept. Electrical Eng.,  \\ 
75013 Paris,& 32000 Haifa  \\ 
France& Israel\\ 
ustunel@enst.fr & zakai@ee.technion.ac.il
\end{tabular}


\begin{thebibliography}{99} 
\markboth{{\sc Bibliography}}{{\sc Bibliography}} 



\bibitem{Carl}
T. Carleman: ``Zur Theorie der linearen
Integralgleichungen''. Math. Zeit. 9, 196-217, 1921.

\bibitem{Del-M}
C. Dellacherie and P. A. Meyer: {\sl Probabilit\'es et Potentiel,
 Ch. I \`a IV}. Paris, Hermann, 1975.


\bibitem{D-S}
N. Dunford and J.T. Schwartz:
{\sl Linear Operators, Vol. 2}, New York, Interscience, 1967.
 
\bibitem{F-P} 
D. Feyel and A. de La Pradelle: ``Capacit\'es gaussiennes''. Annales 
de l'Institut Fourier, {\bf 41}, f. 1, 49-76, 1991. 
 

 
 
\bibitem{FUZ} 
D. Feyel, A.S. \"Ust\"unel and M. Zakai:
``Realization of Positive Random Variables via Absolutely Continuous
Transformations of Measure on Wiener Space''. Probability Surveys,Vol.
3, (electronic) p.170-205, 2006.
 
\bibitem{GIR} 
I.V. Girsanov: ``On transforming a certain class of stochastic
processes by absolutely continuous substitutions of measures''. Theory
of Probability and Appl. {\bf 5}, p. 285-301, 1960.

\bibitem{Gross}
L. Gross: ``Logarithmic Sobolev inequalities''.  Amer. J. Math. {\bf
  97}, no. 4, p.1061--1083, 1975. 

\bibitem{I-W}
N. Ikeda and S. Watanabe:
{\sl Stochastic Differential Equations and Diffusion Processes}. North
Holland, Amsterdam (Kodansha Ltd., Tokyo), 1981.

\bibitem{Kus}
S. Kusuoka: ``The nonlinear transformation of Gaussian measure on
Banach space and its absolute continuity''. J. Fac. Sci., Tokyo Univ.,
Sect. 1 A, {\bf 29}, 567-590, 1982.
 
\bibitem{Mal}
P. Malliavin:
{\sl Stochastic Analysis}. Springer, 1997.

\bibitem{Orey}
S. Orey: ``Radon-Nikodym derivatives of probability measures:
Martingale methods''. Dept. of the Foundations of Mathematical
Sciences, Tokyo University of Education, 1974.

\bibitem{ASU} 
A. S. \"Ust\"unel: 
{\sl Introduction to Analysis on Wiener Space}. 
Lecture Notes in Math. Vol. 1610. Springer, 1995. 

\bibitem{ASU-1} 
A. S. \"Ust\"unel: {\sl Analysis on Wiener Space and
  Applications}. Electronic text at the site
http://www.finance-research.net/.

\bibitem{ASU-2}
A. S. \"Ust\"unel: ``Damped logarithmic Sobolev inequality on the
Wiener space''. Stochastic Analysis and Related Topics VII. The
Silivri Workshop. Progress in Probability, Vol.{\bf 48},
245-249. Birkh\"auser, 2001.

\bibitem{FILT}
A. S. \"Ust\"unel and M. Zakai: ``The construction of filtrations on abstract Wiener space''.
J. Funct. Anal. {\bf  143}  , p. 10--32, 1997.

\bibitem{BOOK} 
A. S. \"Ust\"unel and M. Zakai: 
{\sl Transformation of Measure on Wiener Space}. 
Springer Verlag, 1999. 

\bibitem{CRAS}
A. S. \"Ust\"unel and M. Zakai: ``The invertibility of adapted
perturbations of identity on the Wiener
space''. C. R. Acad. Sci. Paris, S\'erie I, {\bf 342},  p. 689-692, 2006.
 
\end{thebibliography}
\end{document}